\documentclass[10pt,a4paper]{article}
\usepackage{amsmath}
\usepackage{amssymb}
\usepackage{amsthm}
\usepackage[all,2cell,curve,arrow,frame]{xy} \UseAllTwocells
\usepackage{graphicx}
\usepackage{epsfig}
\theoremstyle{definition}

\theoremstyle{remark}

\numberwithin{equation}{section}
\begin{document}
\title{$2$-Categorical Poincar\'{e} Representations and State Sum Applications}
\author{L. Crane \\ {\em Perimeter Institute for Theoretical Physics} \\
{\em Waterloo, ON, Canada} \\ \\
\ M.D. Sheppeard \\
{\em Department of Physics and Astronomy} \\ {\em University of
Canterbury, Christchurch, New Zealand}} \maketitle
\begin{abstract}
This is intended as a self-contained introduction to the
representation theory developed in order to create a Poincar\'{e}
$2$-category state sum model for Quantum Gravity in $4$
dimensions. We review the structure of a new representation
$2$-category appropriate to Lie $2$-group symmetries and discuss
its application to the problem of finding a state sum model for
Quantum Gravity. There is a remarkable richness in its details,
reflecting some desirable characteristics of physical
$4$-dimensionality. We begin with a review of the method of orbits
in Geometric Quantization, as an aid to the intuition that the
geometric picture unfolded here may be seen as a categorification
of this process.
\end{abstract}
\section{Introduction}

There has been much investigation into constrained topological
state sums for Quantum Gravity in four dimensions. The motivation
for these constructions follows largely from the success of
$3$-dimensional models, such as the Turaev-Viro state sum. The
elegant category theoretic formulation of these $3$-dimensional
TFTs has led to a search for even richer $4$-dimensional
analogues.

To date, perhaps the most promising candidates for a finite theory
of quantum gravity are the Lorentzian Barrett-Crane model
\cite{BC} and its q-deformed version \cite{Roche}. In \cite{Crane}
Crane and Yetter propose a new state sum for Lorentzian Quantum
Gravity which utilises a higher algebraic symmetry: the
Poincar\'{e} group as a $2$-category on one object. In order to
attempt to write down the state sum it is first necessary to more
concretely understand the representation category defined in
\cite{Yetter}.

The physical motivation for our construction is the idea that the
fundamental symmetry to use to construct quantum gravity is the
Poincar\'{e} group action, but with the translation subgroup
differentiated from the Lorentz group. Translations could
correspond to lengths, and infinitesimal rotations to bivectors,
the representation category providing a quantization of the whole
system of geometrical quantities. More mathematically expressed,
the Poincar\'{e} group admits a canonical decomposition \[
\mathbf{M}^{4} \rightarrow \mathbf{P} \rightarrow \mathbf{L} \]
where $\mathbf{M}^{4}$ is the additive group of translations in
$\mathbb{R}^{4}$ and $\mathbf{L}$ is the (connected) Lorentz group
$SO(3,1)$; and we need a new notion of representation theory which
respects this decomposition. This information is lost in treating
the Poincar\'{e} group simply as a group.

The way we accomplish this is to regard the Poincar\'{e} group as
naturally a strict $2$-group, which we call $\mathbf{Poinc}$, as
follows. $1$-morphisms are elements of $\mathbf{L}$. $2$-morphisms
$g_{1} \rightarrow g_{2}$ are elements $x \in \mathbf{M}^{4}$ that
map to $g_{1}^{-1} g_{2}$, which in this case takes across all of
$\mathbf{M}^{4}$ when $g_{1} = g_{2}$ and nothing otherwise, as in
the globule
\[ \xymatrix{\ast \rtwocell^{g}_{g}{\hspace{1mm} x} & \ast}
\] The action $\alpha : \mathbf{L}
\times \mathbf{M}^{4} \rightarrow \mathbf{M}^{4}$ appears in the
$2$-category as part of the tensor structure, as does the target
map $t$ taking $\mathbf{M}^{4}$ trivially to the identity in
$\mathbf{L}$.

Higher category theory \cite{Neuchl}\cite{KapVo} has long been
discussed in the context of $4$-dimensional TQFTs. In a more
immediately physical setting it will be interesting to
investigate, for instance, Higher Yang-Mills theories \cite{Baez},
or Higher Lattice Gauge theory \cite{Hendryk}. Unfortunately, the
apparent complexity in the detailed definitions of higher
categories is liable to lead to the misconception that they are ad
hoc and inelegant abstractions. In fact, the construction outlined
here is quite canonical in that it all falls out of one coherent
algebraic structure. We hope to make clear that the new
representation $2$-category is actually a natural object in
differential geometry.

The substitution of a $2$-group for an ordinary group is an
example of the process of categorification \cite{Dolan}, which has
been believed to play an important role in raising the dimension
of categorically constructed theories for some time. Thus the
suggestion that we use a $2$-group as the fundamental symmetry of
a 4D theory is very natural.

In order to understand our new construction, it is useful to first
adopt a slightly novel point of view on the representation theory
of ordinary groups. Any group may be thought of as a category,
with one object $\ast$, with group elements described as
invertible morphisms (arrows). Let us refer to groups so thought
of as $1$-groups. What does the representation theory for
$1$-groups involve? The usual idea of describing a representation
by an action of the group on some vector space is replaced by the
natural categorical analog: a functor from the $1$-group category
into the tensor category of vector spaces. A morphism between such
functors $F$ and $G$ is a  {\em natural transformation}
\cite{MacLane}: for all arrows $g$ there exists $n$ such that
\[ \xymatrix @=8mm {F(\ast) \ar[r]^{F(g)} \ar[d]_{n} & F(\ast) \ar[d]^{n} \\
G(\ast) \ar[r]_{G(g)} & G(\ast) } \] commutes.

The reader should verify that this corresponds to the usual ideas
of group representations and intertwining operators.

Now one would like to construct an appropriate $2$-categorical
analog of the above for $\mathbf{Poinc}$. The most natural
suggestion would be to consider the category of functors between
$\mathbf{Poinc}$ and a suitable $2$-categorical analog of
$\mathbf{Vect}$. Previously studied representation $2$-categories,
such as ${\mathbf 2}$-{\bf Vect} \cite{KapVo} or ${\mathbf
2}$-{\bf Hilb} \cite{BaezTwo}, suffer from the fact that the
corresponding Poincar\'{e} $2$-group representation theory admits
very few representations (as the Lorentz group is not profinite).

In order to remedy this problem, a new $2$-category called
$\mathbf{Meas}$ has been constructed, which we describe below. The
main purpose of this paper is to spell out the concrete structure
of the new $2$-category of representations of the Poincar\'{e}
$2$-group in $\mathbf{Meas}$, introduced in \cite{Yetter},
\cite{BarrMac} and \cite{YetterTwo}. As we shall see, it is very
rich. A large part of the subjects of harmonic analysis, group
representation theory, mathematical physics and low dimensional
topology are combined into a unified algebraic structure in a
peculiarly categorical way. In particular, we observe the
existence of a smooth subcategory closed under tensor product.

Whereas a $1$-category has objects and morphisms, a $2$-category
has objects, $1$-morphisms and $2$-morphisms. For objects $A$ and
$B$ of a $2$-category, Hom$(A,B)$ forms a category. For example,
$\mathbf{Toph}$ is the $2$-category of topological spaces,
homeomorphisms and homotopies.

The category of functors between two $2$-categories is itself a
$2$-category, with functors, weak or strict natural
transformations, and modifications. Thus the new representation
$2$-category has three levels of structure which need to be made
explicit.

Interesting higher categories are not strict but {\em weakened}:
coherence relations for a category may hold only {\em up to} a
$2$-morphism, and coherence relations become higher dimensional
polytopes.  It is precisely the possibility of weakening morphisms
that gives a much richer symmetry out of which one might hope to
build interesting state sums. We shall see below that the weak
natural transformations ($1$-intertwiners) of our $2$-category are
richer geometric objects than the strict ones.

Now let us briefly describe the $2$-category $\mathbf{Meas}$. The
objects of this category \cite{Yetter}\cite{YetterTwo} are the
categories of measurable fields of Hilbert spaces over a given
measure space. $1$-intertwiners are fields of Hilbert spaces on
the product of the domain and range measure spaces, together with
a measure on the product. These should be thought of as matrices
with continuum indices, the measure allowing one to sum up matrix
elements when composing. The $2$-intertwiners are measurable
fields of Hilbert space operators.

A single representation is now a $2$-functor from the $2$-group
into the range $2$-category $\mathbf{Meas}$, and the higher
category of all representations is built out of these functors
with (weakened) pseudo-natural transformations $\{n(X) : F
\rightarrow G \}$ and modifications $\mu : n \Rightarrow m$ as
morphisms. The basic coherence relations for these morphisms are
\[\begin{array}{ccc} {\xymatrix @=20mm {F(X)
\rtwocell^{F(f)}_{F(g)} \ar[d]_{n(X)}
\rtwocell\omit{<8>{\hspace{5mm} n(g)}}
& F(Y) \ar[d]^{n(Y)} \\
G(X) \ar[r]_{G(g)} & G(Y)}} & \simeq & {\xymatrix @=20mm {F(X)
\ar[r]^{F(f)} \ar[d]_{n(X)} \rtwocell\omit{<8>{\hspace{5mm} n(f)}}
 & F(Y) \ar[d]^{n(Y)} \\
G(X) \rtwocell^{G(f)}_{G(g)} & G(Y)}} \end{array} \] and
\[\begin{array}{ccc} {\xymatrix @=20mm {
F(X) \rtwocell^{F(f)}_{F(g)} \rtwocell\omit{<8>{\hspace{5mm}
n(g)}} \dtwocell<-5>^{\hspace{30mm} n(X)}_{m(X)
\hspace{30mm}}{\vspace{4mm} \mu(X)}
& F(Y) \ar[d]^{n(Y)}  \\
G(X) \ar[r]_{G(g)} & G(Y) }} & = & {\xymatrix @=20mm { F(X)
\ar[r]^{F(f)} \ar[d]^{m(X)} \rtwocell\omit{<3>{\hspace{5mm} m(f)}}
& F(Y) \dtwocell<-5>^{\hspace{30mm} n(Y)}_{m(Y)
\hspace{30mm}}{\vspace{4mm} \mu(Y)} \\ G(X)
\rtwocell^{G(f)}_{G(g)} & G(Y)  }} \end{array} \]

By taking the direct integral of the Hilbert spaces in a
$1$-intertwiner of $\mathbf{Meas}$, a representation of the
ordinary Poincar\'{e} group is obtained. However, irreducibles do
not go to irreducibles, so the structure of the representation
$2$-category is new. Despite some technical resemblance it does
not by any means reduce to the famous Wigner classification.

We begin by reviewing Kirillov's {\em method of orbits} for
$1$-groups, as it clarifies somewhat the structure of our new
representation category, which may be thought of as a
categorification of this process.  The analogy between the
structure of the new category and the method of orbits is quite
striking. In place of coadjoint orbits one finds two levels of
orbits: orbits in Minkowski space and fibrations of Lorentz orbits
over these. The tensor product in our new category is a geometric
operation in two steps: the lower level is very close to the
operation of Kirillov. Higher dimensional versions of Geometric
Quantization should make explicit what one means by the {\em
quantum geometry} of $2$-categorical state sums \cite{BaezBarr}.

A note on notation: Objects and functors are denoted by capital
Latin letters, $1$-morphisms by lower-case $f: A \rightarrow B$,
and $2$-morphisms by, for example, $\mu: f \Rightarrow g$.
Monoidal categories have identities $I$ with respect to the
product. The identity arrow of an object $A$ is $1_{A}$.

\section{Geometric Quantization}

The simplest physical example of Geometric Quantization is the
quantization of phase space $\simeq \mathbb{R}^{2n}$ in field
theory. This manifold may be endowed with a symplectic structure,
and all symplectic manifolds look locally like a patch of this
space with coordinates
$(x_{1},x_{2},\cdots,x_{n},p_{1},p_{2},\cdots,p_{n})$. In this
section we briefly outline the general construction.

Let $\mathcal{G}$ denote the Lie algebra of a semi-simple Lie
group $\mathbf{G}$. For our purposes, $\mathbf{G}$ is
$SL(2,\mathbb{C})$, which is locally isomorphic to the Lorentz
group $\mathbf{L}$. An {\em orbit} of the $\mathbf{G}$-action on a
space $M$ has the form of a symmetric space
\[ \mathcal{O} = G / G_{x} \] for $G_{x}$ the isotropy at a point
$x \in M$.

The coadjoint orbits of $\mathbf{G}$ \cite{Kirillov}\cite{Orbits}
have a natural symplectic structure under the action of
$\mathbf{G}$, and those orbits satisfying an integrality condition
induce unitary irreducible representations of $\mathbf{G}$. In
fact, the unitary irreps appearing correspond to the decomposition
of the regular representation.

Functions on $\mathcal{G}^{\ast}$, corresponding to elements of
the Lie algebra, satisfy a Kirillov-Poisson bracket dual to the
Lie algebra structure \cite{Kirillov}. For a basis $X_{i}$ of
$\mathcal{G}$ and structure constants $c^{k}_{ij}$ this bracket is
given by
\begin{equation}
 [f_{1},f_{2}] = \sum_{i,j,k} c^{k}_{ij} X_{k}
\frac{\partial f_{1}}{\partial X_{i}} \frac{\partial
f_{2}}{\partial X_{j}}
\end{equation}
For vector fields on $\mathcal{G}^{\ast}$
\begin{equation}
 v_{i} = \sum_{j,k} c^{k}_{ij} X_{k} \frac{\partial}{\partial X_{j}}
\end{equation}
there is a canonical symplectic structure on each coadjoint orbit
given by
\begin{equation} \omega (v_{i},v_{j}) = \sum_{k} c^{k}_{ij} X_{k} \end{equation}

Generically, full quantization is possible if there exists a
connection $\nabla$ on a line bundle over the orbit \cite{Orbits}
such that the curvature satisfies
\begin{equation} F(\nabla) = 2 \pi i \omega \end{equation}
In other words, the orbit is {\em integral}. This means that the
$1$-dimensional represention of a certain Lie algebra $\mathcal{H}
\subset \mathcal{G}$ may be extended to a unitary representation
of the corresponding group $\mathbf{H}$, and it is from this
representation that a unitary irrep of $\mathbf{G}$ is induced.
This induction step is functorial. That is, we are hoping to view
the representation $2$-category in its geometric guise entirely
within the categorical formalism.

The integrality corresponds to the fact that $\omega \in
H^{2}(M,\mathbb{Z})$, suggesting perhaps that an appropriate
$2$-categorical analogue of quantization would involve higher
cohomological conditions, such as in the theory of gerbes. This
question will not be addressed here.

Explicitly, for the case of $\mathcal{G} = sl(2,\mathbb{C})$ the
orbits on $\mathcal{G}^{\ast}$ are derived as follows. Let
$a,b,c,d$ denote complex variables such that $ad - bc = 1$. Under
the identification $\mathcal{G} \sim \mathcal{G}^{\ast}$ the
coadjoint orbits of $SL(2,\mathbb{C})$ are given by the action of
\[ \left( \begin{array}{ccc}
d^{2} & cd & c^{2} \\
2bd & ad+bc & 2ac\\
b^{2} & ab & a^{2}\end{array} \right) \] on ${\mathbb{C}}^{3}$.
After a suitable change of coordinates, in $\mathbb{R}^{6}$ these
orbits are given by the zero orbit and intersections
\begin{eqnarray} x_{0}^{2} + x_{1}^{2} + x_{2}^{2} - y_{0}^{2} - y_{1}^{2}
- y_{2}^{2} = n^{2} - \rho^{2} \\
x_{0} y_{0} + x_{1} y_{1} + x_{2} y_{2} = n \rho
\end{eqnarray}
for $n, \rho \in \mathbb{R}$. Clearly these curves foliate
$\mathbb{R}^{6}$. Thus we have a remarkably easy classification of
the equivalence classes of unitary representations $\pi_{n \rho}$
for $SL(2,\mathbb{C})$ \cite{Pucz}.

Integrality demands that $n \in \frac{1}{2} \mathbb{Z}$, as for
$SU(2)$, which has spherical orbits in $su(2)^{\ast} \simeq
\mathbb{R}^{3}$. (In the simpler case of $SU(2)$, the integral
orbits correspond to the well known integral levels of total
quantized angular momentum, making rigorous the physicist's
intuition about adding quantum angular momenta by adding vectors
with uncertain direction).

The irreps appearing in the decomposition of a tensor product of
irreps can be recovered geometrically from the sum of orbits
\[ \mathcal{O}_{1} + \mathcal{O}_{2} \equiv
\{ x_{1} + x_{2} : x_{1} \in \mathcal{O}_{1} , x_{2} \in
\mathcal{O}_{2} \}\]

In the case of $SU(2)$, the decomposition $\mathcal{H}_{j} \otimes
\mathcal{H}_{l} = \bigoplus^{|j+l|}_{|j-l|} \mathcal{H}_{i}$ into
spherical shells follows from the range of the norm of the sum of
two vectors of lengths $j$ and $l$. For $SL(2,\mathbb{C})$ we
obtain
\begin{equation} \pi_{n_{1} \rho_{1}} \otimes \pi_{n_{2} \rho_{2}}
= \bigoplus_{m+n_{1} +n_{2} \in \mathbb{Z}} \int^{\oplus} \pi_{m
\rho} d \rho
\end{equation}

This so-called {\em method of orbits} \cite{Kirillov} is seen to
describe the structure of the tensor category of representations
of a Lie group, as a geometric category whose tensor product is a
geometric operation.

The usual unitary representations of the Poincar\'{e} group can be
similarly described by the orbits in Minkowski space along with a
representation of the stabiliser glued to each point of the orbit.

\section{${\mathbf 2}$-Representations}

We now outline the structure of the representation $2$-category,
noting its form as a two stage Geometric Quantization. Its
elements are
\begin{itemize}
\item {\bf Objects: } $2$-functors, labelled by a module object
and an action upon it by $\mathbf{L}$ and $\mathbf{M}^{4}$.
\item {\bf $1$-morphisms: } pseudo-natural transformations
\item {\bf $2$-morphisms: } modifications
\end{itemize}

The general construction, for any $2$-group, was introduced in
\cite{Yetter} and \cite{YetterTwo}. Let $\mathbf{Meas}(X)$ denote
the $1$-category whose objects are measurable fields of Hilbert
spaces $\{ \mathcal{H}_{x} \}$ indexed by the Borel space $X$. The
$1$-morphisms are all bounded fields of bounded operators between
fields of Hilbert spaces. This is well defined if we say a field
of bounded operators is bounded when $x \mapsto \| \phi_{x} \|$,
for $\phi_{x} \in \mathcal{B}(\mathcal{H}_{x},\mathcal{K}_{x})$,
is a bounded real function. The spaces $X$ may be thought of as
continuous generalizations of the discretely labelled categories
$\mathbf{Vect}$-$\mathbf{n}$ underlying
$\mathbf{2}$-$\mathbf{Vect}$.

A single representation is a $2$-functor $\mathcal{R}:
\mathbf{Poinc} \rightarrow \mathbf{Meas}$ from the $2$-group to
the $2$-category of all such categories $\mathbf{Meas}(X)$. The
unique object is mapped to some specific measure space $X$, and
the morphisms, thought of as a category, get mapped to objects and
morphisms in the hom category $\mathbf{Meas}(X,X)$.

More concretely, each element $g$ of $\mathbf{L}$ is sent to a
field of Hilbert spaces on $X \times X$, which must obey the group
law under convolution, and in particular be invertible. In order
for a field of Hilbert spaces to be invertible, each Hilbert space
must be one dimensional, and there can only be one non-zero
Hilbert space in each horizontal or vertical line in $X \times X$.
Thus, the representation of the $1$-intertwiners of
$\mathbf{Poinc}$ means that we are given a measurable action of
$\mathbf{L}$ on the measure space $X$.

Now let us consider the images of the $2$-intertwiners of
$\mathbf{Poinc}$. Each vector of $\mathbf{M}^{4}$ is assigned a
linear map on each Hilbert space in the field corresponding to
each $g \in \mathbf{L}$. These must satisfy the (additive) group
law of the vector space $\mathbf{M}^{4}$. This means that each
point in the graph of the action of each $g \in \mathbf{L}$ is
assigned a character, which may be identified with a point in
$\mathbf{M}^{4}$. The group laws of $\mathbf{Poinc}$ now imply
that the characters are determined by the characters on the graph
of the identity functor on $X$, and that they are equivariant with
respect to the action of $\mathbf{L}$ on $X$.

This translates into the following:

\bigskip
\noindent {\bf Proposition 3.1 } {\em Objects of the $2$-category
of representations of $\mathbf{Poinc}$ in $\mathbf{Meas}$
correspond to measure spaces on which $\mathbf{L}$ acts
measurably, provided with equivariant maps to $\mathbf{M}^{4}$}
\bigskip

For example, let $X$ be an hyperboloid orbit of the Lorentz group
in Minkowski space, such as $\mathcal{O}_{\rho} \equiv \{t^{2} -
x_{1}^{2} - x_{2}^{2} - x_{3}^{2} = \rho \}$. The Lorentz group
acts on $\mathcal{H}_{x}$ by translation on $x \in
\mathcal{O}_{\rho}$, and the characters $\chi \in
\widehat{\mathbf{M}^{4}}$ act as scalar multipliers on the
coordinates of the point.

Now let us classify the irreducible representations. Any orbit of
$\mathbf{L}$ in a representation space $X$ is a subrepresentation.
Orbits of a group in an action correspond to quotients of the
group by the stabilizer subgroup of a point. In order for the
orbit to admit an equivariant map to $\mathbf{M}^{4}$, the
stabilizer subgroup must be contained in the stabilizer of some
point in $\mathbf{M}^{4}$. The image of the quotient under the
equivariant map is then some orbit of $\mathbf{L}$ in
$\mathbf{M}^{4}$.

A simple set of elementary irreducible representations are given
by the orbits $\mathcal{O} \hookrightarrow \mathbf{M}^{4}$, with
the trivial representation attached, on which $\mathbf{L}$ acts
transitively. But these are not the only irreducibles. One must
include multiple copies of the same orbit of the form $M
\rightarrow^{\pi} \mathcal{O}$, where the fibre is a symmetric
space of the stabilizer for the orbit. These copies are permuted
under the action of $\mathbf{L}$.

In other words, irreps correspond to orbits $M$ of
$\mathbf{L}$-actions $\alpha$ which are equivariant fiberings
\[ \xymatrix{M \ar[d]_{\pi} \ar[r]^{\alpha}
& M \ar[d]_{\pi} \\
\mathcal{O}_{\rho} \ar[r] & \mathcal{O}_{\rho} } \]

Because these representations are in $1:1$ correspondence with
subgroups $H \subseteq G_{\rho}$ for each Minkowski orbit type, we
obtain a rather large class of irreducibles, which contains

\begin{enumerate}
\item {\bf Elementary Irreps: }
correspond to the case where the stabilizer of $M$ equals the
stabilizer of the orbit $\mathcal{O}_{\rho}$ of $\mathbf{M}^{4}$
over which it fibers. Denoted by $E_{\rho}$.
\item {\bf Lie Irreps: }
occur where the subgroup of the stabilizer of the orbit
$\mathcal{O}_{\rho}$ is a connected Lie group. For example, if
$\mathcal{O}_{\rho}$ is a spacelike hyperboloid, its stabilizer is
$SU(2)$. One could either take $S^{1}$ for the subgroup, in which
case $M$ is an $S^{2}$ bundle over the hyperboloid, or let the
stabilizer be trivial, in which case the space $M$ is a copy of
the group $\mathbf{L}$, which can be written as an $S^{3}$ bundle
over the hyperboloid. This maximal irrep, which can occur over any
orbit in $\mathbf{M}^{4}$, we denote by $L_{\rho}$.
\item {\bf Crystallographic Irreps: }
occur if the stabilizer of $M$ contains a discrete subgroup of the
stabilizer of the orbit in $\mathbf{M}^{4}$. We obtain an irrep
whose fiber over a point in $\mathbf{M}^{4}$ is a manifold whose
fundamental group includes the given discrete group.
\item {\bf Non-Hausdorff Irreps: }
occur when we choose as stabilizer a non-Lie subgroup of
$\mathbf{L}$ to produce a total space which is not Hausdorff.
\end{enumerate}

There is a subset of irreps, namely the first three cases, which
inherit a natural smooth structure. This subcategory will close
under tensor product and direct integral by a smooth index space.
The study of the non-Hausdorff irreps will require very different
mathematical tools.

\bigskip
\noindent {\bf Proposition 3.2 } {\em For all orbits of
$\mathbf{L}$ in $\mathbf{M}^{4}$, the 4 cases above exhaust the
irreducible objects of $\mathbf{Rep(Poinc)}$}
\bigskip

\subsection{Tensor Products of Objects}

The tensor product of two objects $M_{1}$ and $M_{2}$ in
$\mathbf{Meas}$ corresponds to the cartesian product of the
underlying measure spaces. The $2$-group has a natural action on
the tensor product of two representations, where the functors
corresponding to elements of $\mathbf{L}$ act in both variables at
once, while the actions of the $2$-morphisms of $\mathbf{Poinc}$
on the spaces over a point in $M_{1} \times M_{2}$ is the tensor
product of the respective actions. Note that the elementary zero
orbit $E_{0}$ acts as an identity element.

Now the tensor product of two characters of the abelian group
$\mathbf{M}^{4}$ just corresponds to vector addition. Thus the
tensor product of two objects of $\mathbf{Rep(Poinc)}$ corresponds
to taking the sum of their two projections as subsets of the
vector space $\mathbf{M}^{4}$ and projecting the product space
$M_{1} \times M_{2}$ into $\mathbf{M}^{4}$ by sending each point
to the vector sum of the projections of its two coordinates. Note
the analogy with the method of orbits.

The decomposition of tensor products into direct integrals of
irreps is easy enough to work out explicitly in the simpler cases,
uncovering an interesting picture. Whereas products of basis
elements in an algebra are determined by structure coefficients,
and tensor products in a category based in ${\mathbf 2}$-{\bf
Vect} are defined by structure vector spaces, the tensor product
of irreps in $\mathbf{Rep(Poinc)}$ decomposes in terms of
structure spaces. These structure spaces are measure spaces in
general, but smooth manifolds for the most natural cases. Observe
that in this $2$-categorical theory, representations are therefore
{\em not} linear at all levels. Linearity is only required at the
top level, where we wish to define trace operators.

We denote this type of direct integral
\begin{equation} \int^{\uplus} (M,F) \end{equation}
where $F$ is the fibre of $M$.

Let us illustrate this by means of a simple but important example:
the case of the tensor product of two elementary irreps,
corresponding to orbits $\mathcal{O}_{1}$ and $\mathcal{O}_{2}$ in
$\mathbf{M}^{4}$. The fibers over each point of each orbit are
single points. Therefore, the fiber of the tensor product over any
point of $\mathbf{M}^{4}$ is the set of all ways to decompose the
point as a sum of two points, one in $\mathcal{O}_{1}$ and one in
$\mathcal{O}_{2}$. The decomposition into irreps is accomplished
by first decomposing into orbits in $\mathbf{M}^{4}$, then
decomposing the fiber over a generic point into orbits of the
stabilizer of the orbit.

Consider the case where both $\mathcal{O}_{1}$ and
$\mathcal{O}_{2}$ are spacelike hyperboloids with radii $\rho_{1}$
and $\rho_{2}$. The set-theoretic sum $\mathcal{O}_{1} +
\mathcal{O}_{2}$ is the union of all spacelike hyperboloids with
radius $\rho \geq \rho_{1} + \rho_{2}$. The fiber over each point
is the set of all timelike triangles in $\mathbf{M}^{4}$ with side
lengths $(\rho_{1}, \rho_{2}, \rho)$. Generically this set forms a
$2$-sphere, which is a single orbit of $SU(2)$. Thus the tensor
product of two spacelike irreps is a direct integral of Lie irreps
with stabilizer $S^{1}$
\begin{equation} E_{\rho_{1}} \otimes E_{\rho_{2}} =
(\int^{\uplus}_{\rho > \rho_{1} + \rho_{2}}
 ( M_{\rho} , S^{2} )) + E_{\rho_{1} + \rho_{2}}
\end{equation}
The second term corresponds to the collinear $(\rho_{1}, \rho_{2},
\rho)$ case.

Extending this procedure to the tensor product of several
elementary irreps allows us to compute the tensor product of the
various Lie irreps.

The triple tensor product contains copies of the maximal irreps
$L_{\rho}$ of allowable orbits because it contains the space of
quadrilaterals with appropriate edge lengths (figure 1). In other
words, a generic quadrilateral in $\mathbf{M}^{4}$ admits no
isotropy subgroup. In this decomposition there is not simply one
maximal irrep, but rather a family, indexed by the space of shapes
of such quadrilaterals in $\mathbf{M}^{4}$. Copies of
$(M_{\rho},S^{2})$ arise for planar quadrilaterals, of isotropy
group $SO(2)$.

\begin{figure}
\centering \epsfbox{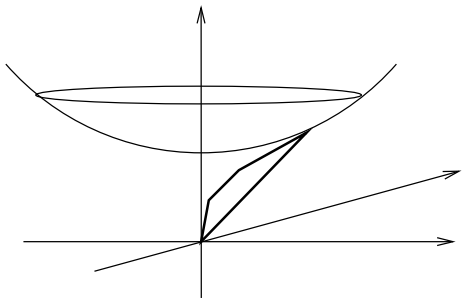} \caption{quadrilateral in
$\mathbf{M}^{4}$}
\end{figure}

Observe that the structure spaces for such cases are smooth, and
the irreps that appear in the decomposition of tensor products of
fundamental irreps are also smooth. This suggests that it is
possible to restrict to a smooth subcategory, which might be
expedient in applications. Moreover, in studying deformations, or
cohomology, of our category, in the smooth case everything will
reduce to a tractable study of suitable forms on products of shape
spaces.

The tensor products of all Lie irreps will be worked out in detail
in later papers. The problem for crystallographic irreps will need
slightly different methods. In some cases, this family is quite
large, i.e. corresponds to Fuchsian groups. It is remarkable that
an algebraic structure combining them is defined canonically. The
non-Hausdorff irreps will probably be harder to describe
explicitly. They do not appear in the tensor products of the
smooth representations.

\section{${\mathbf 1}$-Intertwiners}

A {\em strong} $1$-intertwiner between two objects in
$\mathbf{Rep(Poinc)}$ is given by a field of Hilbert spaces on the
product of their underlying measure spaces, which is invariant
under the product group operation and confined to ordered pairs
which fiber over the same point $x$ in $\mathbf{M}^{4}$. That is,
\begin{equation} \mathcal{H}_{x} \mapsto \int dx \hspace{2mm} \mathcal{H}_{x}
\otimes \mathcal{K}_{(x,y)} d\mu_{y}(x)
\end{equation} which is a translation of
\begin{equation}  \begin{array}{ccc} {\xymatrix @=20mm
{X \rtwocell^{\mathcal{R}_{1}(f)}_{\mathcal{R}_{1}(g)} \ar[d]
\rtwocell\omit{<8>{\hspace{5mm} n(g)}}
& X \ar[d] \\
Y \ar[r]_{\mathcal{R}_{2}(g)} & Y}} & = & {\xymatrix @=20mm {X
\ar[r]^{\mathcal{R}_{1}(f)} \ar[d] \rtwocell\omit{<8>{\hspace{5mm}
n(f)}} & X \ar[d] \\
Y \rtwocell^{\mathcal{R}_{2}(f)}_{\mathcal{R}_{2}(g)} & Y}}
\end{array}
\end{equation}

A {\em weak} $1$-intertwiner has the same structure, together with
a linear action of the isotropy group of $x$ on the Hilbert space
assigned to it.

Thus, in order to obtain an explicit description of the
irreducible $1$-intertwiners between two representations we
decompose the cartesian product of the underlying measure spaces
into group orbits and describe the orbits and the orbit space. In
the weak case we include irreducible representations of the
isotropy group as well. The general (weak) $1$-intertwiner is then
a direct integral of any measurable family of points in the orbit
space.

In other words, for irreps of the form
\[ \xymatrix @=5mm {& F_{i} \ar[r] & M_{i} \ar[d]_{\pi} \\
&& \mathcal{O}_{\rho} } \] which are symmetric spaces for
$\mathbf{L}$, $1$-intertwiners $f: M_{1} \rightarrow M_{2}$ must
be equivariant with respect to the action of $\mathbf{L}$, i.e.
$f$ is concentrated on the fibres $F_{1}$ and $F_{2}$. Let \[
F_{1} \sim \frac{G_{\mathcal{O}_{\rho}}}{G_{1}} \hspace{1.5cm}
F_{2} \sim \frac{G_{\mathcal{O}_{\rho}}}{G_{2}} \] Then
intertwiners are built on bridges
\begin{equation} \xymatrix @=5mm { & \mathcal{I} \ar[dl] \ar[dr] & \\
F_{1} && F_{2}} \end{equation} for \[ \mathcal{I} \sim
\frac{G_{\mathcal{O}_{\rho}}}{G_{1} \cap G_{1}} \]

In particular, note that the space of irreducible $1$-intertwiners
between two smooth irreps is smooth. This further justifies the
idea of a smooth subcategory.

In the simplest case of $1$-intertwiners from an elementary irrep
$E_{\rho}$ to itself, we obtain only the constant functions on the
orbit, since coherence implies the pointwise condition
\begin{equation}
n(g_{1},x) n(g_{2},g_{1}(x)) = n(g_{1}g_{2},x)
\end{equation}
and we are free to choose a basis at each $x \in E_{\rho}$ under
which $n$ must be invariant, forcing $n \equiv 1$.

This is very suggestive with regards to the problem of
construction of state sum models: the weak $1$-intertwiners are
the space of integrable spinors at a point, i.e. the
representations of the isotropy group at the point. For example,
for $(M_{\rho},S^{2}) \rightarrow (M_{\rho},S^{2})$ we obtain
fields of measurable subsets of $S^{2}$.

The composition of $1$-intertwiners requires that the convolution
of two orbits be decomposed into orbits. For tensor product, the
product of orbits must be decomposed into orbits. In general, the
tensor product of two $1$-intertwiners is defined weakly by a
choice of $2$-intertwiner \cite{Crans}
\begin{equation} \xymatrix @=20mm {{X_{1} \otimes Y_{1}} \ar[r]^{X_{1} \otimes g}
\ar[d]_{f \otimes Y_{1}}
& {X_{1} \otimes Y_{2}} \ar[d]^{f \otimes Y_{2}} \\
{X_{2} \otimes Y_{1}} \ar[r]_{X_{2} \otimes g} \utwocell\omit{<8>
f \otimes g} & {X_{2} \otimes Y_{2}}}
\end{equation}

\subsection{Internal Hom Functor and Duality}

The new $2$-category has, modulo some subtleties about measures,
an internal Hom functor and a duality. These satisfy the usual
relationship with the tensor product \begin{equation}
\mathrm{Hom}(A,B \otimes C) = \mathrm{Hom}(A \otimes B^{\ast}, C)
\end{equation} The duality comes from reflection about the origin
in $\mathbf{M}^{4}$. The internal Hom functor comes from
identifying equivariant Hilbert fields with Hilbert fields over
the space of orbits under $\mathbf{L}$.

In the case of two irreducibles fibering over the same orbit in
$\mathbf{M}^{4}$, the orbit space should be thought of as fibering
over the origin in $\mathbf{M}^{4}$. This makes perfect sense in
terms of the above formula, if we let $C$ be the trivial
representation $E_{0}$.

The subtlety is that a $1$-intertwiner requires a choice of
measure. If we only want to work in the smooth subcategory, we can
restrict to the usual Lebesgue measures on the appropriate spaces
and the above relation becomes rigorous.

These structures are necessary for the detailed definition of a
state sum model.

\section{${\mathbf 2}$-Intertwiners}

A $2$-morphism in $\mathbf{Meas}$ is defined as a measurable field
of operators on the cartesian product space on which the two
$1$-morphisms are defined, from one field of Hilbert spaces to the
other. In order for a $2$-morphism between two $1$-intertwiners to
be a modification of functors, and hence a $2$-intertwiner in
$\mathbf{Rep(Poinc)}$, it must intertwine the characters $\chi \in
\widehat{\mathbf{M}^{4}}$. Thus, a $2$-intertwiner between
irreducible $1$-intertwiners is given by a measurable scalar
function on the orbit in the cartesian product, where it is
defined. This all takes place within a fiber over
$\mathbf{M}^{4}$. The horizontal and vertical compositions
\[ \xymatrix{X \rtwocell^{\mathcal{R}(f)}_{\mathcal{R}(f)}
& X \rtwocell^{\mathcal{R}(g)}_{\mathcal{R}(g)} & X} \hspace{20mm}
\xymatrix{X \ruppertwocell^{\mathcal{R}(f)} \ar[r]
\rlowertwocell_{\mathcal{R}(f)} & X}
\] are defined by convolutions. For example, for two constant
function intertwiners on $\mathcal{O}_{\rho}$, convolution is the
usual convolution of functions.

The conditions on $2$-intertwiners from tensor products \cite{GPS}
are
\begin{equation} \begin{array}{ccc} {\xymatrix @=20mm {
{X_{1} \otimes Y_{1}} \rtwocell^{X_{1} \otimes f}{\hspace{6mm}
X_{1} \otimes \beta} \rtwocell\omit{<8>{\hspace{5mm} h \otimes k}}
\dtwocell<-5>_{g \otimes Y_{1} \hspace{30mm}}{\vspace{4mm} \alpha
\otimes Y_{1}}
& {X_{1} \otimes Y_{2}} \ar[d] \\
{X_{2} \otimes Y_{1}} \ar[r] & {X_{2} \otimes Y_{2}} }} & = &
{\xymatrix @=20mm { {X_{1} \otimes Y_{1}} \ar[r] \ar[d]
\rtwocell\omit{<6>{\hspace{5mm} g \otimes f}} & {X_{1} \otimes
Y_{2}}
\dtwocell<-5>^{\hspace{30mm} h \otimes Y_{2}}{\alpha \otimes Y_{2}} \\
{X_{2} \otimes Y_{1}} \rtwocell_{X_{2} \otimes k}{\hspace{6mm}
X_{2} \otimes \beta} & {X_{2} \otimes Y_{2}} }}
\end{array}
\end{equation}
and
\begin{equation} \begin{array}{ccc} {\xymatrix @=8mm
{{X_{1} \otimes Y_{1}} \ar[r] \ar[d]
\xtwocell[1,1]{}\omit{\hspace{6mm} f_{1} \otimes g_{1}} & {X_{1}
\otimes Y_{2}} \ar[r] \ar[d] \xtwocell[1,1]{}\omit{\hspace{6mm}
f_{1} \otimes g_{2}}
& {X_{1} \otimes Y_{3}} \ar[d] \\
{X_{2} \otimes Y_{1}} \ar[r] \ar[d]
\xtwocell[1,1]{}\omit{\hspace{6mm} f_{2} \otimes g_{1}} & {X_{2}
\otimes Y_{2}} \ar[r] \ar[d] \xtwocell[1,1]{}\omit{\hspace{6mm}
f_{2} \otimes g_{2}} &
{X_{2} \otimes Y_{3}} \ar[d] \\
{X_{3} \otimes Y_{1}} \ar[r] & {X_{3} \otimes Y_{2}} \ar[r] &
{X_{3} \otimes Y_{3}}} } & = & {\xymatrix @=15mm {{X_{1} \otimes
Y_{1}} \ar[r] \ar[d] \xtwocell[1,1]{}\omit{\omit (f_{2}f_{1})
\otimes (g_{2}g_{1})}
& {X_{1} \otimes Y_{3}} \ar[d] \\
{X_{3} \otimes Y_{1}} \ar[r] & {X_{3} \otimes Y_{3}}}}
\end{array}
\end{equation}

\section{Topological Invariants and Quantum Gravity}

It is expedient to bear in mind that the representations form a
$3$-category on one object \cite{GPS} (the monoidal structure
allows us to do this). Thus edges in a state sum get mapped to
$2$-functors, as seems natural.

The original motivation for this research \cite{Crane} was to
construct a higher categorical version of the Lorentzian
Barrett-Crane model \cite{BC}\cite{PReview} whose state sum, based
on the dual $2$-complex, is
\begin{equation} \mathcal{Z}_{BC} = \sum_{color}^{\int} \prod_{f \in \Delta_{2}}
(\rho_{f}^{2} + n_{f}^{2}) \prod_{e \in \Delta_{1}}
\mathcal{A}_{e} \prod_{v \in \Delta_{0}} (15j)_{v}
\end{equation}
where the sum actually contains an integral over the continuous
parameter $\rho$.

We see that the irreps of $\mathbf{Rep(Poinc)}$ contain a choice
of orbit in $\mathbf{M}^{4}$, ie. a radius. This radius can be
interpreted as a length, either positive or negative, on an edge
in a state sum. The function spaces which appear in the structure
of $\mathbf{Rep(Poinc)}$ can be decomposed into representations of
$\mathbf{L}$, providing a connection to the Barrett-Crane model.
The question of appropriate constraints on the topological theory
remains to be studied, but here we outline the basic structure of
the new state sum.

The topological Poincar\'{e} state sum is constructed by labelling
edges with representations from $\mathbf{Meas}$, faces with
$1$-intertwiners and tetrahedra with $2$-intertwiners. That is, a
functor from a suitable category of triangulated $4$-manifolds, or
PL pseudomanifolds, to the category $\mathbf{Meas}$, a braided
version of which will become a $4$-category with one object and
one $1$-morphism.

The amplitude of a $4$-simplex is calculated by tracing over the
five tetrahedra of its boundary via one of two possible paths
\begin{eqnarray}
+: \alpha_{2} \otimes \alpha_{4} \rightarrow \alpha_{1} \otimes
\alpha_{3} \otimes \alpha_{5} \\
-: \alpha_{1} \otimes \alpha_{3} \otimes \alpha_{5} \rightarrow
\alpha_{2} \otimes \alpha_{4}
\end{eqnarray}
The condition of {\em sphericity} \cite{Mackaay} on
$\mathbf{Rep(Poinc)}$ says that these two maps are equal.

\bigskip
\noindent {\bf Proposition 6.1 } {\em $\mathbf{Rep(Poinc)}$ is a
spherical $2$-category}
\bigskip

We now examine the structure on a tetrahedron in one simple case.
First, for edges of a closed triangle
\[ \xymatrix{& \ast \ar[dr]^{E_{2}} & \\ \ast
\rrlowertwocell<0>_{E_{3}}{<-2>f} \ar[ur]^{E_{1}} && \ast} \]
labelled by three elementary spacelike irreps, such that
$(\rho_{1},\rho_{2},\rho_{3})$ forms an allowable triangle, a
$1$-intertwiner
\[ f: E_{1} \otimes E_{2} \rightarrow E_{3} \] is a map $S^{2} \rightarrow
\bullet$. Replacing $E_{3}$ with $(M_{3},S^{2})$ we obtain all
equivariant maps $S^{2} \rightarrow S^{2}$.

The full labelling of a tetrahedron amounts to the choice of two
$1$-intertwiners, each a composition of two faces of the
tetrahedron, and a choice of $2$-intertwiner between them. Tracing
over tetrahedra gives $4$-simplex labels, which we call $5j$
symbols.

More specifically, in the case of timelike edges we get the set of
$2$-intertwiners from one copy of the product of two $(M_{\rho} ,
S^{2})$ to another, corresponding to the two halves of the
tetrahedron. This function space decomposes into a tensor product
of functions on the orbits in $\mathbf{M}^{4}$ with the functions
on the $S^{2}$ fibers. Other choices about the causal structure of
the tetrahedron would give analogous results, based on orbit
spaces for the stabilizer subgroups $SO(2,1)$, or $E(2)$ in the
null case. If we decompose these into harmonics, we will get a
combination of relativistic balanced spin nets, as in \cite{BC},
with ordinary spin nets. The exact form of the connection will
contain the constraints for the new state sum, which has yet to be
figured out.

The categorical state sum seems to be leading us towards a
plausible way of approaching the quantization of the geometry of a
triangulated $4$-manifold: quantize the edges and configuration
spaces of the triangles independently, then constrain them. We
would not have hit on this procedure without the $2$-category as a
guiding framework. It remains to be seen how well it will work.

Thus a $2$-categorical topological state sum for the smooth
subcategory, formally at least, looks like
\begin{equation} \mathcal{Z} = \mathcal{N} \sum^{\int}_{color}
\prod_{e \in \Delta_{1}} \rho_{e} \prod_{f \in \Delta_{2}}
\mathcal{A}_{f} \prod_{t \in \Delta_{3}} \mathcal{A}_{t} \prod_{s
\in \Delta_{4}} (5j)_{s}
\end{equation}
where $\mathcal{A}_{f}$ follows from the decomposition of function
spaces into representations of $\mathbf{L}$ and their orbit
Casimirs, in analogy with the Lorentzian BC model.

\section{Conclusions}

In the course of working on the structure of
$\mathbf{Rep(Poinc)}$, we were surprised to discover how rich it
is. A plethora of new possibilities for the construction of state
sum models appears to arise. It is natural to speculate on the
possible significance of the Lie and crystallographic irreps for
features of unification such as phase transitions, or even the
insertion of matter.

The application of the new structure to quantum geometry has not
been worked out yet, but looks promising. The fibered spaces in
the objects we put on parts of triangulations have natural
interpretations as staged quantizations of geometrical variables.
Problems of regularization remain to be considered, as do the
appropriate constraints for quantum gravity.

We have focused on the Poincar\'{e} $2$-group for reasons of
physical interest. In fact, our construction has an analog for any
choice of a semisimple Lie group and a representation of it
\cite{Baez}. This means that a whole new chapter of representation
theory is opened up, closely allied to a new family of quantum
geometries.

The discovery of quantum groups has played a very important role
recently in a number of areas of mathematics and physics. Quantum
groups can be approached by applying deformation theory to the
category of representations of a Lie group. Now we have an
interesting $2$-category of representations of a $2$-group. The
techniques of deformation theory have recently been extended to
tensor $2$-categories in \cite{elgueta}. It would be interesting
to study the deformations of $\mathbf{Rep(Poinc)}$ with, for
example, topological invariants in mind.

In \cite{DaySt} it is also observed that for quantum group
$1$-categories it is really the comodules that should be regarded
as the representation theory. The $2$-categorical analog is
considered, and this would be an interesting avenue for further
investigation. The utility of perverse sheaves in understanding
quantum groups presumably has a higher dimensional analog, which
might be approached from the aforementioned deformation theory.

We believe that when a higher braiding is added, to obtain a
$4$-category with two singular levels, that it will be possible to
begin considering unified models wherein the higher category
theory dictates a subtle interplay between massive and spacetime
degrees of freedom.

In its modern form, the functorality of Geometric Quantization
promises a rigorous approach to the quantum geometry of
$4$-dimensional state sums. It encompasses, for instance, theorems
\cite{Moment} on the commutativity of quantization and symplectic
reduction, and on non-Abelian localization, such as formulae
derived by Witten \cite{Witten} from the path integral of
$2$-dimensional Yang-Mills theory.

This paper exhibits but one step on the road to an explicit
construction of a Poincar\'{e} state sum, which is a promising
candidate for a physical theory of quantum gravity. It should be
pointed out that further categorification is possible in defining
$4$-dimensional spin foam models. Higher category theory is as yet
poorly understood, and in particular there are various definitions
of weak $n$-categories, the relationship between them still
somewhat mysterious.

\bigskip
\noindent {\bf Acknowledgements} \\
Both authors are thankful for discussions with H. Pfeiffer, J.
Baez, D. Yetter, L. Freidel, L. Smolin, F. Girelli and others at
the Perimeter Institute, Waterloo, where this work was carried
out. The second author is partly supported by the Marsden Fund of
the Royal Society of New Zealand.
\bibliographystyle{unsrt}
\bibliography{Stash}
\end{document}